\documentclass[10pt,twoside,reqno]{article}
\usepackage{amsmath, amssymb}
\usepackage{epsfig}
\usepackage{cite}
\usepackage{url}

\usepackage{anysize}
\marginsize{2.0cm}{2.0cm}{2.0cm}{2.0cm}
\newtheorem{thm}{Theorem}[section]

\numberwithin{equation}{section}
\newcommand{\R}{\mathbb{R}}

\newcommand{\twofractals}[4]{%
\begin{center}
\begin{minipage}[t]{0.45\textwidth}
\begin{centering}
\includegraphics[width=0.98\textwidth]{#1}
\end{centering}
\centering{\small  #2}
\end{minipage}\hfill%
\begin{minipage}[t]{0.45\textwidth}
\begin{centering}
\includegraphics[width=0.98\textwidth]{#3}
\end{centering}
\centering{\small #4}
\end{minipage}\hfill%
\end{center}%
}
\newcommand{\threefractalss}[6]{%
\begin{center}
\begin{minipage}[t]{0.31\textwidth}
\begin{centering}
\includegraphics[width=0.98\textwidth]{#1}
\end{centering}
\centering{\small #2}
\end{minipage}\hfill%
\begin{minipage}[t]{0.31\textwidth}
\begin{centering}
\includegraphics[width=0.98\textwidth]{#3}
\end{centering}
\centering{\small #4}
\end{minipage}\hfill%
\begin{minipage}[t]{0.31\textwidth}
\begin{centering}
\includegraphics[width=0.98\textwidth]{#5}
\end{centering}
\centering{\small #6}
\end{minipage}
\end{center}%
}
\begin{document}
\setcounter{page}{1}
\begin{center}
\vspace {0.4cm}
\vspace{0.4cm}
{\large{\bf A new tool to study real dynamics: \emph{The Convergence Plane}}}\\
\vspace{0.4cm}
\'Angel Alberto Magre\~n\'an\\
Universidad de La Rioja\\
Departamento de Matem\'aticas y Computaci\'on\\
26002 Logro\~no, La Rioja, Spain\\
alberto.magrenan@gmail.com\\
\end{center}

\begin{abstract}
In this paper, the author presents a new tool, called \emph{The Convergence Plane}, that allows to study the real dynamics of iterative methods whose iterations depends on one parameter in an easy and compact way. This tool can be used, inter alia, to find the elements of a family that have good convergence properties and discard the bad ones or to see how the basins of attraction changes along the elements of the family. To show the applicability of the tool an example of the dynamics of the Damped Newton's method applied to a cubic polynomial is presented.
\vspace{0.4cm}
\\
\end{abstract}
\section{Introduction and Motivation}\label{Intro}
The main aim of the author in this paper is to present a new tool which can make easier the study of the real dynamics of families of iterative methods which depends on a certain parameter or even the study of an iterative method applied to a uniparametric family of polynomials. This tool can be modified in order to extend, amongst other ones, to methods which needs two approximations as for example, secant-type methods, modified Newton's method, etc.

The dynamics of iterative methods used for solving nonlinear equations in complex plane has been studied recently by many authors \cite{2,3,4,5,6,12,13,14,15}. There exists a belief that real dynamics is included in the complex dynamics but this is not true at all. For example in the real dynamics one can proof the monotone convergence which does not exist in the complex plane, there exists also asymptotes in the real dynamics but in the complex one that concept has no sense or the point $z=\infty$ in the complex plane can be studied as another point but in the real line it is not possible. As a consequence real dynamics is not contained in the complex dynamics and both must be studied separately. Taking into account that distinction, many authors \cite{1,7,9,11,13} have begun to study it since few years ago.

If one focus the attention on families of iterative methods studied in the complex plane, parameter spaces have given rise to methods whose dynamics are not well-known. These parameter spaces consist on studying the orbits of the free critical points associating each point of the plane with a complex value of the parameter. Several authors \cite{4,5,7,13} have studied really interesting dynamical and parameter planes in which they have found some anomalies such as convergence to $n$-cycles, convergence to $\infty$, or even chaotical behavior. In the real line, there exists tools such as Feigenbaum diagrams or Lyapunov exponents that allow us to study what happens with a concrete point, but it is really hard to study each point in a separate way. This is the main motivation of the author to present the new tool called \emph{The Convergence Plane}.

\emph{The Convergence Plane} is obtained by associating each point of the  plane with a value of the starting point and a value of the parameter. That is, \emph{The Convergence Plane} is based on taking the vertical axis as the value of the parameter and the horizontal axis as the starting point, so every point in the plane represents an initial estimation and a member of the family. If one draws a straight horizontal line in a concrete value of the parameter, the dynamical behavior, for that value, for every starting point is on that line. On the other hand, if the straight line is vertical, the dynamical behavior for that starting point and every value of the parameter is on that line, this is the information that gives Feigenbaum diagrams or Lyapunov exponents, so both tools are included.

The rest of the paper is organized as follows: in Section \ref{S2} the Algorithm of \emph{The Convergence Plane} is shown and in Section \ref{S3} an example of \emph{The Convergence Plane} associated to the Damped Newton's method applied to the polynomial $p_-(x)=x^3-x$ is provided in order to validate the tool. Finally, the conclusion are shown in the concluding Section \ref{S4}.

\section{Algorithm of \emph{The Convergence Plane}}\label{S2}

As it is said in the introduction each point of the plane corresponds to a starting point and a value of the parameter, in other words the pair $(x_0, \lambda_0)$ represents that the study is developed using $x_0$ as the starting point and $\lambda_0$ as the value of the parameter. The algorithm of this new tool is the following one:
\begin{itemize}
  \item The fixed points of the method must be computed.
  \item Then, a color is assigned to each fixed point.
  \item Moreover, the region in which one wants to study the family, $D$, the maximum number of iterations, $M$,  and the tolerance, $\varepsilon$, are prefixed.
  \item Then, a grid of $d\times d$ points in $D$ of the initial points and values of the parameter must be chosen.
  \item If after $M$ iterations of the family with $\lambda_0$ as the value of the parameter, the point $x_0$ does not converge to any of the fixed points that point must be black.
  \item If one is interested on representing the convergence to any cycle or other behavior such as convergence to extraneous fixed points, divergence etc., there exists the possibility of assigning a color to that behaviors too. \end{itemize}

Once \emph{The Convergence Plane} has been computed it is easy to distinguish the pairs $(x_0,\lambda_0)$ for which the element of the family is convergent to any of the roots using $x_0$ as a starting point. So this tool provides a global vision about what points converges and shows what are the best choices of the parameters to ensure the greatest basin of attraction. Moreover, it can be used also as a tool that show how the basins of attraction of a family changes with the value of the parameter or even to study the convergence of methods which depend on 2 points such as, for example, secant-type methods.
\subsection{Mathematica program}
In this moment, the author is going to explain how Figure \ref{Fig1} and Figure \ref{Fig2} of this paper were generated. To do this, it is shown the Mathematica programs that has been used  which are a modification of the ones that appears in \cite{16}. In concrete, in the example the region will be $[-2,2]\times [2,2.26]$ (the season of taking that values is in \cite{13} and will appear in Section 3). Moreover, a number of iterations $M=1000$,  a tolerance of $\varepsilon=10^{-6}$ and a grid of $1024\times1024$ have been taken.

First of all, the function, the roots and the procedure that identifies the root to which converge the iterations have been defined as
\begin{verbatim}
f[x_]:=x^3-x;

rootf[1]=-1;
rootf[2]=0;
rootf[3]=1;

rootPosition = Compile[{{z,_Real}},
    Which[
    Abs[z - rootf[1]] <  10.0^(-6), 1,
	    Abs[z - rootf[2]] <  10.0^(-6), 2,
	    Abs[z - rootf[3]] <  10.0^(-6), 3,
	    Abs[z ] >  10.0^(3), 11,
	    True, 0],
    {{rootf[_],_Real}}
    ]
\end{verbatim}
Then, the iteration method is the following
\begin{verbatim}
iterNewtonLamda=Compile[{{x,_Real},{k,_Real}},x-k*f[x]/f'[x]]
\end{verbatim}
The algorithm used to show if the iteration of \emph{iterMethod} of a point converges to a root or a cycle is the following
\begin{verbatim}
 iterColorAlgorithm[iterMethod_,xx_,yy_,lim_] :=
    Block[{z,z2,kk,ct,r}, z =xx; kk=yy;
    ct = 0;r = rootPosition[z];
    While[(r==0) && (ct < lim),++ct; z = iterMethod[z,kk];
    r = rootPosition[z];];
    If[r==0,z2= iterMethod[z,kk];z2= iterMethod[z2,kk];
        If[Abs[z2-z]<10^(-6),r=4,z2=iterMethod[z2,kk];
          If[Abs[z2-z]<10^(-6),r=5,z2=iterMethod[z2,kk];
            If[Abs[z2-z]<10^(-6),r=6,z2=iterMethod[z2,kk];
              If[Abs[z2-z]<10^(-6),r=7,z2=iterMethod[z2,kk];
                If[Abs[z2-z]<10^(-6),r=8,z2=iterMethod[z2,kk];
                  If[Abs[z2-z]<10^(-6),r=9,z2=iterMethod[z2,kk];
                    If[Abs[z2-z]<10^(-6),r=10,z2=iterMethod[z2,kk];
                      ]]]]]]]]
      If[Head[r]==Which,r =0]; (* "Which" unevaluated *)
      Return[N[r+ct/(lim+0.001)]]
    ]
\end{verbatim}
The palette of colors used is defined as
\begin{verbatim}
ConvergenceColor[p_] :=
    Switch[IntegerPart[11p],
                       	         11, CMYKColor[0.0,0.0,0.0,0.0],(*White*)
                                 10, CMYKColor[0.0,1.0,1.0,0.8],(*Dark red*)
                                  9, CMYKColor[0.0,0.0,1.0,0.5],(*Dark yellow*)
                                  8, CMYKColor[1.0,0.0,1.0,0.5],(*Dark Green*)
                                  7, CMYKColor[1.0,0.0,0.0,0.5],(*Dark blue*)
      	                           6, CMYKColor[0.0,0.5,1.0,0.0],(*Orange*)
      	                           5, CMYKColor[1.0,0.0,1.0,0.0],(*Green*)
                                  4, CMYKColor[0.0,1.0,1.0,0.0],(*Red*)
      	                           3, CMYKColor[0.0,0.0,1.0,0.0],(*Yellow*)
                                  2, CMYKColor[0.0,1.0,0.0,0.0],(*Magenta*)
                                  1, CMYKColor[1.0,0.0,0.0,0.0],(*Cyan*)
                                  0, CMYKColor[0.0,0.0,0.0,1.0](*Black*)
              ]
\end{verbatim}
The function used by the author to plot the convergence plane is the following
\begin{verbatim}
plotConvergencePlane[iterMethod_,points_] :=
        DensityPlot[iterColorAlgorithm[iterMethod,x,k,limIterations],
	        	  {x, xxMin, xxMax}, {k, kkMin, kkMax},
	        	  PlotRange->{0,11}, PlotPoints->points, Mesh->False,
	        	  ColorFunction->ConvergenceColor]
\end{verbatim}
Then, a graphic is obtained in this way (notice that we avoid overflow and underflow errors and other errors by means of using the instruction \emph{Off}):
\begin{verbatim}
  numberPoints = 1024; limIterations = 4000;
  xxMin = -2.0; xxMax = 2; kkMin = 0.0; kkMax = 2.6;
  Off[General::ovfl];             Off[General::unfl];           Off[Infinity::indet]
  Off[CompiledFunction::cccx];    Off[CompiledFunction::cfn];   Off[CompiledFunction::cfcx];
  Off[CompiledFunction::cfex];    Off[CompiledFunction::crcx];  Off[CompiledFunction::cfse];
  Off[CompiledFunction::ilsm];    Off[CompiledFunction::cfsa];
  plotConvergencePlane[iterNewtonLamda, numberPoints]
  \end{verbatim}

\section{Example: Damped Newton's method applied to the polynomial $p_-(x)=x^3-x$}\label{S3}

To show the goodness of \emph{The Convergence Plane} we are going to study the dynamics of the Damped Newton's method which has the following form
\begin{equation}\label{1}
N_{\lambda,p}(x)=x-\lambda\dfrac{p(x)}{p'(x)},
\end{equation}
where $p(x)$ is a polynomial with real coefficients and we take $\lambda\in\R$. In concrete, in \cite{13}, the author has made an extensive and deep study about the real dynamics of the damped Newton's method applied to polynomials of degrees 2, 3, 4 and 5. In order to proof that the tool works properly we are going to apply the iterative method \eqref{1} to find the real roots of a cubic polynomial. The Scaling Theorem \cite{1,2,3,7,11,13} allows up to suitable change of coordinates, to reduce the study of the dynamics of iterations of general maps, to the study of specific families of iterations of simpler maps. Specifically, the study of cubic polynomials reduce to the study of $p_0(x)=x^3$, $p_+(x)=x^3+x$ $p_-(x)=x^3-x$ and the uniparametric family $p_\gamma(x)=x^3+\gamma x+1$. In \cite{13} the author  shows the real dynamics of the Damped Newton's method applied to $p_-(x)=x^3-x$ are not easy and that there exists some cycles and chaotic behavior of the iterations of some points by means of using Lyapunov exponents,  Feigenbaum diagrams and analytical techniques. We are going to use the new tool \emph{The Convergence Plane} to study the dynamics of the Damped Newton's method applied to $p_-(x)$. We denote the three roots of $p_-(x)$ as $r_1=-1$, $r_2=0$ and $r_3=1$.

\begin{figure}[h!]
\twofractals{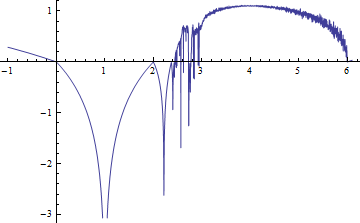}{Lyapunov Exponents}%
{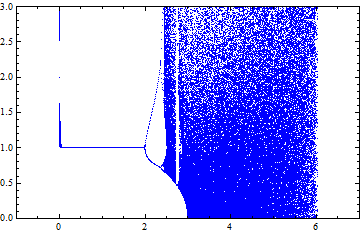}{Feigenbaum diagram}%
\caption{Lyapunov exponents and Feigenbaum diagram associated to the iteration of Damped Newton's method applied to the polynomial $p_-(x)=x^3-x$, on the interval $\lambda \in (-1,7)$ taking $x_0=2.0$.}
\label{FFF1}
\end{figure}

Before, applying the new tool we have to see what information give to us the Lyapunov exponents and Feigenbaum diagrams. These two tools gives information only about the point which is being iterated. In concrete, the Lyapunov exponent is defined, for an orbit $\{x_1,\ldots,x_n,\ldots\}$ as
$$
h(x_1)=\displaystyle{\lim_{n\to\infty}}\frac{1}{n}\left({\log{|f'(x_1)|}+\cdots+\log{|f'(x_n)|}}\right).
$$
The applicability of this tool resides in the following result.
\begin{thm}\cite{9}
An  orbit $\{x_1,\ x_2,\ \ldots,\ x_n\}$ is chaotic if and only if the following conditions hold:
\begin{itemize}
\item The orbit is not asymptotically periodic.
\item $h(x_1)>0$.
\end{itemize}
\end{thm}
As a consequence, this tool shows which orbits are chaotical or not. On the other hand, the Feigenbaum diagrams shows if the orbit of a point converges to a cycle, to a point, if it is chaotic or it diverges. In Figure \ref{FFF1} it is shown the results obtained using the both two tools.

\begin{figure}[h!]
\twofractals{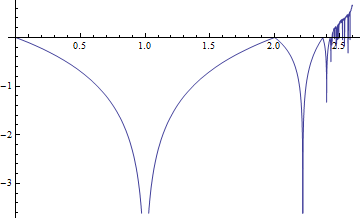}{Lyapunov Exponents}%
{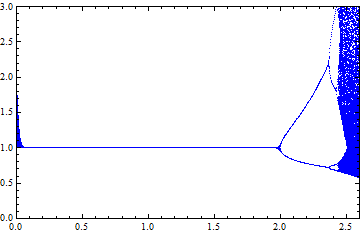}{Feigenbaum diagram}%
\caption{Lyapunov exponents and Feigenbaum diagram associated to the iteration of Damped Newton's method applied to the polynomial $p_-(x)=x^3-x$ on the interval $\lambda \in (0,2.6)$ taking $x_0=2.0$.}
\label{FF1}
\end{figure}

Taking into account that a brief study has been made in \cite{13} and the idea of not making this paper very long we are going to focus on the interval $\lambda\in(0,2.6)$ which is sufficient to show the powerful of the new tool. In Figure \ref{FF1} we can see both tools centered on that interval. And we distinguish three different zones.

\begin{figure}[h!]
\threefractalss{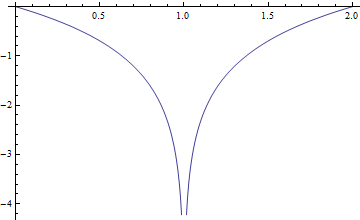}{$\lambda\in(0,2]$}%
{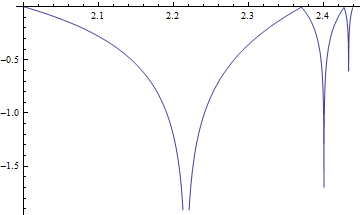}{$\lambda\in(2,2.439]$}%
{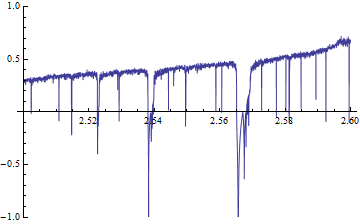}{$\lambda\in(2.5,2.6]$}%
\threefractalss{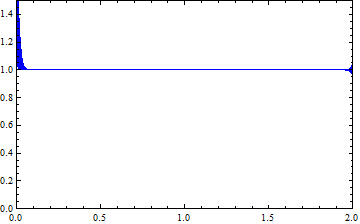}{$\lambda\in(0,2]$}%
{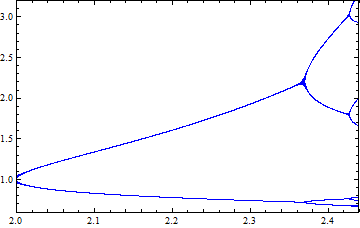}{$\lambda\in(2,2.439]$}%
{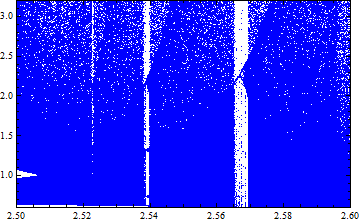}{$\lambda\in(2.5,2.6]$}%
\caption{Feigenbaum diagrams and Lyapunov exponents associated to the iteration of Damped Newton's method applied to the polynomial $p_-(x)=x^3-x$ on different intervals of $\lambda$ taking $x_0=2.0$.}
\label{FF2}
\end{figure}

There exists 3 clear zones of different behavior. In the first one which corresponds with the interval $\lambda \in (0,2)$ (see left side of the Figure \ref{FF2}), we see the zone in which the iterations converge to the fixed points of the Damped Newton's method or equivalently to the roots of the polynomial $p_-(x)$. The second zone, shown in the center of Figure \ref{FF2} is where it appears cycles of different orders, in concrete cycles of order $2$, $4$ and $8$. The third zone, shown in the right side of Figure \ref{FF2} is where we found chaotical behavior or convergence to different cycles of periods different to $2^n$.

\begin{center}
\begin{figure}[htpb!]
\begin{center}
  \includegraphics[width=0.7\textwidth]{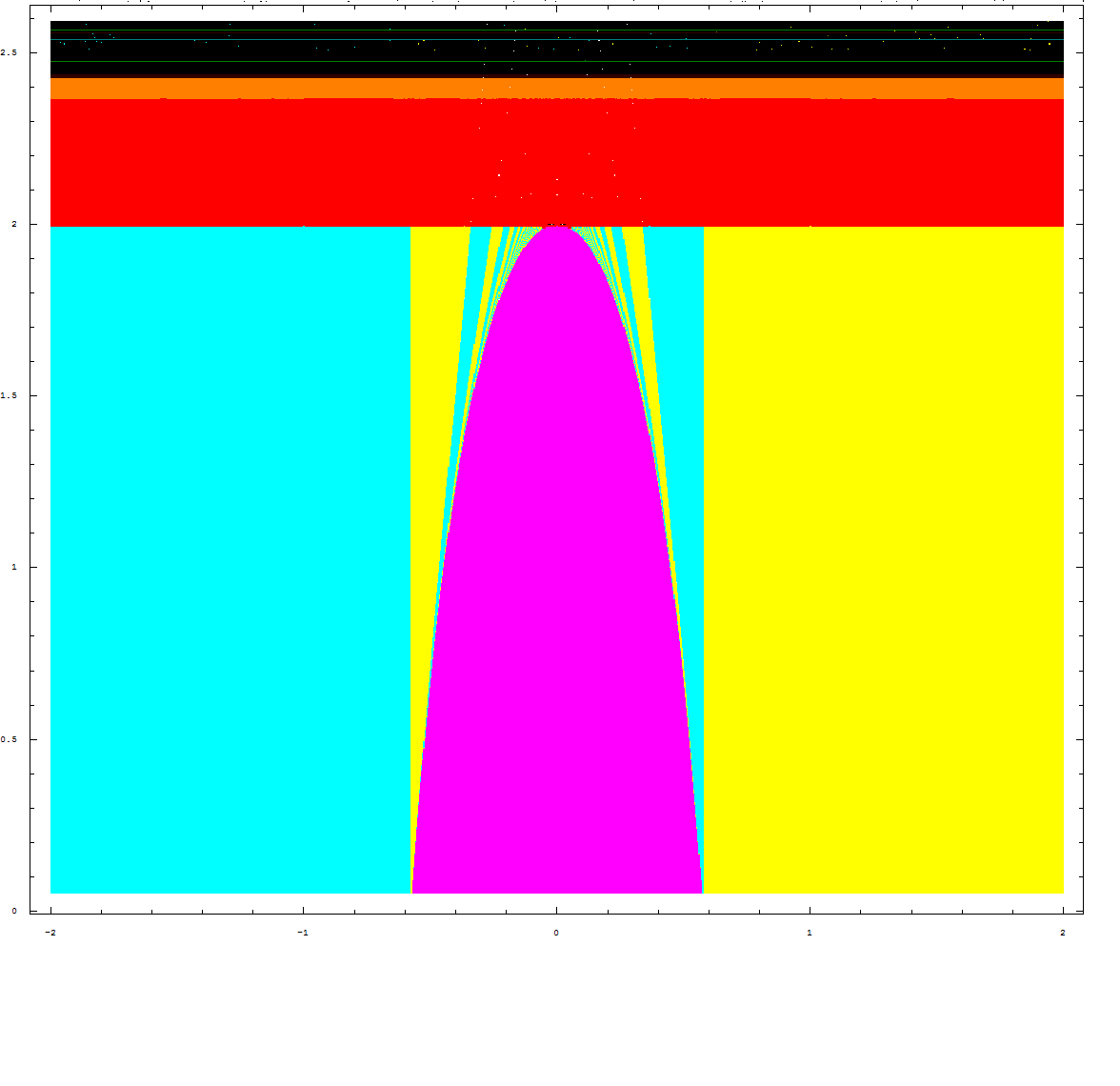}\\
  \caption{\emph{The Convergence Plane} associated to the Damped Newton's method applied to the polynomial $p_-(x)=x^3-x$ on the region $(x,\lambda)\in[-2,2]\times[0,2.6]$.}\label{Fig1}
\end{center}
\end{figure}
\end{center}

Now if we use \emph{The Convergence Plane} in the region $[-2,2]\times[0,2.6]$ (see Figure \ref{Fig1}) the goodness of this tool is going to be proof. In this case, we use the program \emph{Mathematica 5.0} (the code appears in Section 2.1) as in \cite{10,16}, with tolerance $\varepsilon = 10^{-6}$, a maximum of 4000 iterations and the following palette of colors:
\begin{itemize}
  \item Cian, if the iterations converge to root $r_1$, magenta, if the iterations converge to root $r_2$ and yellow, if the iterations converge to root $r_3$.
  \item Red, if the iterations converge to a $2$-cycle, green, if converge to a $3$-cycle, orange, if converge to a $4$-cycle, dark blue, if converge to a $5$-cycle, dark green, if to a $6$-cycle, dark yellow, if converge to a $7$-cycle and dark red if the iterations converge to a $8$-cycle.
  \item White, if the iterations diverge to $\infty$.
  \item Black, in other case.
\end{itemize}

\begin{figure}[h!]
\threefractalss{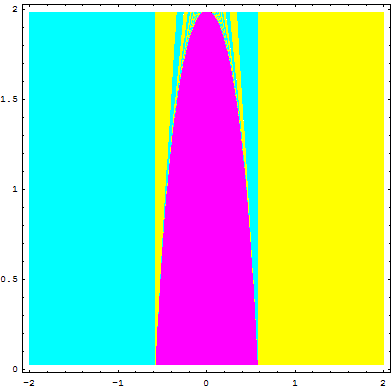}{$\lambda\in(0,2]$}%
{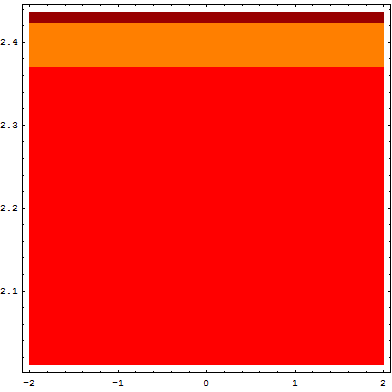}{$\lambda\in(2,2.439]$}%
{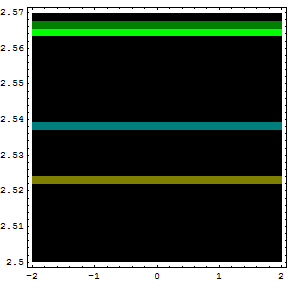}{$\lambda\in(2.5,2.6]$}%
\caption{The 3 zones with different behavior of the iterations of Damped Newton's method for different values of the parameter $\lambda$.}
\label{Fig2}
\end{figure}

Again we distinguish 3 clear zones with different dynamical behavior, which corresponds with the 3 zones using Feiganbaum diagrams and Lyapunov exponents. The first one, shown in the left hand of Figure \ref{Fig2} corresponds to the interval $\lambda \in (0,2]$ and it is the zone in which every point (except the poles and its preimages), converges to any of the roots of the polynomial or equivalently to any fixed point of the Damped Newton's method. Additionally, in this zone the tool gives the idea of how the basins of attraction changes with the value of $\lambda$. For example, the basin of $r_2=0$ is getting lower when the value of $\lambda$ is closer to two, an the basin of the other two roots are getting bigger as the parameter is closer to $2$. Furthermore, it it shown that the Julia set gets more intricate when the value of the parameter increases until $2$. The second zone, see the center of Figure \ref{Fig2}, corresponds with the zone in which there exists convergence to cycles of orders $2$, $4$ and $8$. And the third zone, shown in the right hand of the Figure \ref{Fig2}, the zone in which there exists chaotical behavior and convergence to cycles of order different to $2^n$. In concrete in the interval $\lambda \in (2.5,2.6)$ appears cycles of order $5$, $6$ and $7$.

Summarizing, the author has shown the helpfulness of this new tool using it to study the behavior of the Damped Newton's method applied to a cubic polynomial. The conclusions drawn this study are: the Damped Newton's method is a good method as a root-finding algorithm if $\lambda\in(0,2]$, if $\lambda\not\in(0,2]$, the iteration can converge to cycles, diverge to infinity, etc.; the basin of attraction of $r_2=0$ increases when the damping factor decreases, and the other basin increases with the parameter; the Julia set associated with the Damped Newton's method applied to $p_-(x)=x^3-x$ is more intricate as the damping factor increases to $2$. All these conclusions coincide with the ones given in \cite{13}.

\section{Conclusions}\label{S4}

In the present paper, the author has presented a new tool that allows to study the convergence of a family of iterative methods, taking into account every initial point. \emph{The Convergence Plane}, includes the information given by tools such as Lyapunov exponents and Feigenbaum diagrams, indeed, the new tool provides more information because it considers each initial point of the real line and every value of the parameter. Moreover, this technique can be used to choose the best value of the parameter that ensure the convergence zone as large as possible and shows how the basins of attraction changes with the value of the parameter. On the other hand, this technique can be modified in order to: study methods applied to polynomials which depends on a parameter, to other kind of methods such as two-point methods (for example secant-type methods), methods applied to non-differentiable functions, etc. Therefore, this new tool can be used in a large amount of situations, making the study of real dynamics easier, deeper and in a more compact way.

\end{document}